\documentclass{elsart}
\journal{Advances in Applied Mathematics}
\usepackage{amsmath,amssymb,verbatim,amscd,graphicx}

\newcommand\C{\mathbb C}

\newcommand\PP{\mathbb P}

\newcommand\diag{\operatorname{diag}}

\newcommand\rank{\operatorname{rank}}

\renewcommand\O{\mathcal O}

\newcommand\aideal{\mathfrak a}
\newcommand\bid{\mathfrak b}
\newcommand\st{*}

\begin{document}

\begin{frontmatter}



\title{Phylogenetic ideals and varieties for the general Markov
model}

\date{July 18, 2006}


\author{Elizabeth S. Allman},
\ead{e.allman@uaf.edu}
\author{John A. Rhodes}
\ead{j.rhodes@uaf.edu}

\address{Department of Mathematics and Statistics\\University of
Alaska Fairbanks\\Fairbanks, AK 99775}

\begin{abstract}
The general Markov model of the evolution of biological sequences
along a tree leads to a parameterization of an algebraic variety.
Understanding this variety and the polynomials, called phylogenetic
invariants, which vanish on it, is a problem within the broader area
of Algebraic Statistics. For an arbitrary trivalent tree, we
determine the full ideal of invariants for the 2-state model,
establishing a conjecture of Pachter-Sturmfels. For the
$\kappa$-state model, we reduce the problem of determining a
defining set of polynomials to that of determining a defining set
for a 3-leaf tree. Along the way, we prove several new cases of a
conjecture of Garcia-Stillman-Sturmfels on certain statistical
models on star trees, and reduce their conjecture to a family of
subcases.

\end{abstract}


\end{frontmatter}

\section{Introduction}

An important problem arising in modern biology is that of
sequence-based \emph{phylogenetic inference.} Suppose we obtain a
collection of biological sequences, such as genomic DNA, from
currently extant species, or \emph{taxa}. Assuming these sequences
evolved from a common ancestral sequence, how can we infer a tree
that describes their evolutionary descent? The use of algebraic
methods for this problem was first proposed in 1987 in independent
works by Lake \cite{Lake87}, and Cavender and Felsenstein
\cite{CF87}. Recently, Garcia, Stillman, and Sturmfels \cite{GSS}
initiated a more general algebraic study of statistical models, of
which phylogenetic models are a particularly interesting example. In
this new field, Algebraic Statistics, the viewpoints of algebraic
geometry are central to investigations of probabilistic models
arising in applied contexts.

\smallskip

In model-based phylogenetics, evolution is usually assumed to
proceed along a binary tree from an ancestral sequence at the root
of the tree, to sequences found in the taxa, which label the leaves
of the tree. The $\kappa=4$ bases $A,C,G,T$ of which DNA is composed
are viewed as states of random variables. Each site in the sequence
might be assumed to evolve i.i.d., so that different sites can be
viewed as trials of the same process. Probabilities of the various
base substitutions along an edge of the tree can then be given by a
Markov transition matrix along that edge.  Additional biologically
reasonable, or mathematically convenient, assumptions as to the form
of these transition matrices are often imposed. The basic problem is
to assume some model along these lines and use it to infer, from
observations of DNA sequences only at the leaves, a tree topology
that might describe their evolutionary descent.  An excellent
overview of the field of phylogenetics is provided by the recent
volume of Felsenstein \cite{Fel04}.

\smallskip

In the phylogenetics literature, a \emph{phylogenetic invariant} for
a particular model and tree is a polynomial that vanishes on all
joint distributions of bases at the leaves that arise from the
model, regardless of the values of the model parameters. In the
terminology of algebraic geometry, the model and tree imply a
parameterization of a dense subset of a variety, and phylogenetic
invariants are the elements of the prime ideal defining that
variety.

For applications, one might hope that the near-vanishing of
phylogenetic invariants on observed frequencies of bases in DNA data
could be used as a test of model-fit and/or tree topology. Although
this idea remains undeveloped for practical use, phylogenetic
invariants have already provided means for addressing more
theoretical questions in phylogenetics, such as the nature of
maximum likelihood points \cite{CHHP}, and the identifiability of
certain models \cite{ARidtree}.

\medskip

In this paper we investigate the phylogenetic variety for the
general Markov model of base substitution for an arbitrary tree, a
detailed specification of which will be given in the next section.
This model was also the focus of the related investigations
\cite{AR03,ARQuart}.

One main result is the proof of Conjecture 13 of Pachter and
Sturmfels \cite{PS} on the ideal of phylogenetic invariants for the
general Markov model in the case of $\kappa=2$ states: the
invariants arising from all $3\times 3$ minors of `2-dimensional
flattenings' of an array along the edges of a binary $n$-taxon tree
$T$ generate the full ideal.  This is Theorem \ref{thm:PS}, which is
stated more fully in Section \ref{sec:SPconj} and proved in Section
\ref{sec:k2id}.

\begin{figure}[h]
\begin{center}
\includegraphics[height=1in, width=1.35in]{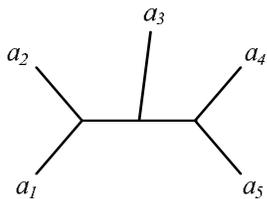}
\end{center}
\caption{A 5-taxon tree}\label{fig:5taxa}
\end{figure}

For an explicit example of this theorem, consider the 5-taxon tree
of Figure \ref{fig:5taxa}.  Then for the 2-state model, denote the
states by 0 and 1. A $2\times 2\times 2\times 2\times 2$ tensor $P$
encodes the probabilities of various states at the leaves, where
$P(i_1,i_2,i_3,i_4,i_5)=p_{i_1i_2i_3i_4i_5}$ is the joint
probability of observing state $i_j$ in the sequence at leaf $a_j$,
$j=1,\dots,5$.  Now $P$ has two natural flattenings according to the
partitions of leaves produced by deleting an internal edge of the
tree. The partitions, or \emph{splits}, are
$\{\{a_1,a_2\},\{a_3,a_4,a_5\}\}$, and
$\{\{a_1,a_2,a_3\},\{a_4,a_5\}\}$, and the corresponding flattenings
are
$$\begin{pmatrix} p_{00000}& p_{00001}& p_{00010}& p_{00011}& p_{00100}&
p_{00101}& p_{00110}& p_{00111} \\ p_{01000}& p_{01001}& p_{01010}&
p_{01011}& p_{01100}& p_{01101}& p_{01110}& p_{01111}\\ p_{10000}&
p_{10001}& p_{10010}& p_{10011}& p_{10100}& p_{10101}& p_{10110}&
p_{10111}\\ p_{11000}& p_{11001}& p_{11010}& p_{11011}& p_{11100}&
p_{11101}& p_{11110}& p_{11111} \end{pmatrix}$$
and
$$\begin{pmatrix}
p_{00000}&p_{00001}&p_{00010}&p_{00011}\\
p_{00100}&p_{00101}&p_{00110}&p_{00111}\\
p_{01000}&p_{01001}&p_{01010}&p_{01011}\\
p_{01100}&p_{01101}&p_{01110}&p_{01111}\\
p_{10000}&p_{10001}&p_{10010}&p_{10011}\\
p_{10100}&p_{10101}&p_{10110}&p_{10111}\\
p_{11000}&p_{11001}&p_{11010}&p_{11011}\\
p_{11100}&p_{11101}&p_{11110}&p_{11111} \end{pmatrix}.  $$ The
theorem states that the $3\times 3$ minors of these two matrices
generate the prime ideal of all phylogenetic invariants for the
2-state general Markov model on this tree. In particular, this ideal
has a natural set of generators that correspond to the splits, and
therefore to specific topological features of the tree.

We note that this theorem provides the first determination of all
phylogenetic invariants for an arbitrary binary tree for any
non-group-based model. Sturmfels and Sullivant
\cite{q-bio.PE/0402015} solved the similar problem for group-based
models, using the Hadamard conjugation
(\cite{Hen89,HenPen89,MR93m:62121,MR1218244}) to recognize the
varieties as toric. While algebraic models intermediate to the
group-based and general ones have been introduced recently
\cite{ARSBD,math.AG/0407033}, our knowledge of them is less
complete.

\smallskip

We also investigate the question of the explicit determination of
the phylogenetic variety and ideal for larger $\kappa$. We show in
Theorem \ref{thm:kany} that if we have a set of polynomials whose
zero set is the variety for the 3-taxon tree, then we can construct
a set of polynomials whose zero set is the variety for any binary
$n$-taxon tree. Similar to the conjecture of \cite{PS}, our
constructions involve `flattenings', though both 2- and
3-dimensional ones are now needed, as might be expected from
\cite{AR03}.  Thus the only remaining obstruction to our
determination of a defining set of polynomials for the phylogenetic
variety for any binary tree and any number of states $\kappa$ is the
determination of a defining set for the 3-taxon tree variety.

In Conjecture \ref{con:phylotree} we suggest that the same
construction yielding set-theoretic defining polynomials for the
variety would yield generators of the full prime ideal vanishing on
the variety, provided we begin with generators of the ideal for the
3-taxon tree. This is the analog for arbitrary $\kappa$ of the
Pachter-Sturmfels conjecture.

Theorem \ref{thm:PS}, Theorem \ref{thm:kany}, and
Conjecture \ref{con:phylotree}, as well as the Sturmfels-Sullivant
group-based result, can all be viewed as statements that
the phylogenetic varieties and ideals arise from the `local structure'
of the tree. Exploiting this observation to provide better ways of
characterizing the statistical support a data set might provide for
specific local tree features would be interesting work for the future.
In particular, invariants might provide a means of characterizing
support for particular splits or tripartitions of the taxa.
\smallskip

Despite our primary focus on phylogenetic models, to prove Theorem
\ref{thm:kany} we must consider certain other statistical models on star
trees. In Section \ref{sec:startrees}, we therefore investigate models
with a $\kappa$-state hidden variable associated to the internal node,
and $l_i$-state observed variables associated to the $n$ leaves. Such
models are of course interesting in applications outside of
phylogenetics, as they are examples of rather common `mixture models' in
statistics. Following \cite{GSS}, they are termed \emph{hidden naive
Bayes models}.

Our work here focuses on such models in the case that for each $i$ the
number of states $l_i$ is at least as large as the number of hidden
states $\kappa$.  Theorems \ref{thm:ststar} and \ref{thm:idstar}
describe how set-theoretic and ideal-theoretic defining sets of the
associated varieties can be deduced from set-theoretic and
ideal-theoretic defining sets of the variety of the related model which
has $\kappa$-state variables on each leaf.

As a consequence of this work on star tree models,
in Corollary \ref{cor:GSS345}
we prove several cases of
Conjecture 21 of \cite{GSS}, on ideal generators for the hidden naive
Bayes model with $\kappa=2$.  While one of these cases, for the 3-leaf
tree, has been recently proved in \cite{lands}, even for that case our
argument is different, and perhaps more direct. Moreover, our work indicates
that establishing the special
cases mentioned in Conjecture \ref{con:GSSsimp} of this paper
is sufficient to prove
the full conjecture of \cite{GSS}.

\medskip

Before obtaining these results, we begin with several background
sections.  In Section \ref{sec:apvar}, we define the phylogenetic
variety for the general Markov model through the natural
parameterization arising from modeling molecular evolution along a tree
$T$ by associating Markov matrices to each edge.  In Section
\ref{sec:reparam} we then give a more convenient parameterization of (a
dense subset of) the cone over the phylogenetic variety, which
associates an arbitrary $\kappa\times \kappa$ matrix to each edge of
$T$, rather than a Markov matrix. Section \ref{sec:SPconj} introduces
flattenings of tensors along edges and vertices of trees, while Section
\ref{sec:alg} develops the relationship of a form of multiplication of
tensors to the varieties under investigation.  Subsequent sections
contain our primary results.

\medskip
Finally, we note that most of the results on phylogenetic trees in
this paper hold not only for binary trees, but also under the weaker
assumption that each vertex have valency at least three. An
important exception is Theorem \ref{thm:PS}, where the binary
assumption is critical to our proof.

\section{Affine and Projective Phylogenetic Varieties}
\label{sec:apvar}

Let $T$ denote an $n$-taxon tree, by which we mean a tree with all
internal vertices unlabeled and of valency at least 3, with $n$
leaves labeled by taxa $a_1, \dots,a_n$. We will sometimes specify
in addition that $T$ is binary ({\it i.e.}, all internal vertices
are trivalent), as this assumption is needed for some of our
results, and is often the case of primary interest in phylogenetics.

Choosing as a root any vertex $r$ of $T$, either internal or a leaf,
denote the rooted tree by $T^r$.  Parameters for the $\kappa$-state
\emph{general Markov model} of sequence evolution on $T^r$ consist
of a root distribution vector
$\boldsymbol{\pi}_r=(\pi_1,\pi_2,\dots,\pi_\kappa)$ with
non-negative entries summing to 1, together with a $\kappa\times
\kappa$ Markov matrix $M_e$, which has non-negative entries with
each row summing to 1, for each of the $2n-3$ edges $e$ of $T^r$
directed away from $r$.

This models the evolution of biological sequences as follows.  The
$\kappa$ states $[\kappa]=\{1,2,\dots,\kappa\}$ correspond to the
alphabet from which sequences are composed.  The root $r$ represents the
most recent common ancestor of the currently extant taxa, and other
internal nodes of the tree represent most recent common ancestors of
those taxa separated from the root by that node.  The root distribution
vector encodes the frequencies $\pi_i$ with which each state $i$
occurs in an ancestral sequence at $r$. The $(i,j)$-entry of a Markov
matrix along a particular edge of $T$ directed away from $r$ is the
conditional probability of state $i$ changing to state $j$ at any
particular site in the sequence during evolution along that edge. Thus
each site in a biological sequence is assumed to evolve independently,
according to the same process (i.i.d.).  Note the biological term
`sequence' as used here implies no mathematical structure other than a
correspondence between sites based on ancestry; except for matching
sites by common ancestry, the ordering of the sites within the sequences
is irrelevant.

Suppose a rooted $n$-taxon tree $T^r$ has $|E|$ edges, so that for a
binary tree $|E|=2n-3$. For the general Markov model of evolution
along $T^r$ the parameter space $S$ can thus be identified with a
subset of $[0,1]^N$, where $N=(\kappa-1)+|E|\kappa(\kappa-1)$.

Furthermore, there is a polynomial map $\phi_r:S \to [0,1]^L$,
$L=\kappa^n$, giving the joint distribution of states in sequences
at the leaves resulting from any parameter choice. We view points in
$ \phi_r(S)$ or $\C^L$ as $\kappa\times\cdots\times\kappa$ tensors,
with the $i$th index referring to the state at leaf $a_i$. Indices
thus typically range through $[\kappa]$, and a fixed ordering of the
taxa is reflected in the ordering of indices of tensors.  Assuming
the model adequately reflects real molecular evolution, from
biological sequence data we can estimate entries of $\phi_r(s)$, but
usually have little or no direct information about the parameters
$s$.

\medskip

The map $\phi_r$ is explicitly given by $\phi_r(s)=P$, where
\begin{equation}
P(i_1,\dots, i_n)=\sum_{(b_v)\in H} \left ( \boldsymbol{\pi}_r(b_r)
\prod_{e}M_e(b_{s(e)},b_{f(e)})\right ),\label{eq:Pdef}
\end{equation}
where the
product is taken over all edges $e$ of $T^r$ directed away from $r$,
edge $e$ has initial vertex $s(e)$ and final vertex $f(e)$ and
associated Markov matrix $M_e$, and the sum is taken over the set $$H=\{
(b_v)_{v\in Vert(T)} ~|~ b_v\in [\kappa] \text{ if $v\ne a_j$},\ b_v=i_j
\text{ if $v=a_j$}\}\subset [\kappa]^{2n-2}.$$ Thus $H$ represents the
set of
all `histories' consistent with the specified states at the leaves.

\smallskip

The map $\phi_r$ can also be defined inductively, using matrix
algebra, by viewing the tree $T^r$ as built up from smaller trees by
the addition of pairs of terminal edges, as we now explain. For this
purpose, we first assume $T$ is binary.

A \emph{cherry} of $T$ is a pair of distinct leaves
$a_{i_1},a_{i_2}$ whose incident edges contain a common (internal)
vertex of $T$. For $n\ge 3$, any binary $n$-taxon tree contains at
least two cherries, and any rooted binary $n$-taxon tree contains at
least one cherry in which neither taxon of the cherry is the root of
the tree.

For $n\ge 3$ let $T^r_n=T^r$ denote a rooted binary $n$-taxon tree
labeled by taxa $a_1,\dots, a_n$. Choose a cherry of $T^r_n$ which
does not contain the root $r$. Let $T^r_{n-1}$ denote the rooted
binary $(n-1)$-taxon tree obtained by deleting the cherry and its
two incident edges from $T_n^r$ and labeling as a new taxon, say
$b$, the (formerly internal) common vertex of the incident edges.

Applying this definition recursively, we obtain from $T^r$ a sequence of
rooted trees $T^r_n,T^r_{n-1},\dots, T^r_2$, which of course may depend
on some arbitrary choices of cherries. We assume such choices have been
made and fixed.

The map $\phi_r$ described above can now be described inductively as
follows:

A rooted $2$-taxon tree has only one edge $e$ directed away from $r$, so
with parameters $\boldsymbol{\pi}_r$ and $M_e$,
$$\phi_r(\boldsymbol{\pi}_r;M_e)=\diag(\boldsymbol{\pi}_r)M_e,$$ where
$\diag(\mathbf v)$ denotes the square matrix with $\mathbf v$ on its
main diagonal and zeros elsewhere.

To define $\phi_r$ for $T^r_m$, direct edges $e$ away from $r$ and
suppose parameters $$s=(\boldsymbol{\pi}_r; \{ M_e \})$$ for $T^r_m$
are given. Then one obtains parameters $\tilde s$ for $T^r_{m-1}$ by
simply discarding from $s$ the two Markov matrices associated to the
edges of $T^r_m$ not appearing in $T^r_{m-1}$. Inductively, we may
assume $\tilde \phi_r:\Tilde S\to [0,1]^{\kappa^{m-1}}$, the map
giving the joint distribution of states at leaves for $T^r_{m-1}$ as
a function of parameters on $T^r_{m-1}$, has been given. For
convenience, we also assume that taxa of $T^r_m$ are $a_1,a_2,\dots,
a_m$ and those of $T^r_{m-1}$ are $a_1,a_2,\dots, a_{m-2},b$, with
the given orderings, and that $e_1$ and $e_2$ are the edges of
$T^r_m$ containing $a_{m-1},a_m$ respectively.

Then $\phi_r(s)=P$, where $P$ is an $m$-dimensional tensor with
2-dimen\-sional slices given by first letting $\Tilde P=\tilde
\phi_r(\tilde s)$, $\mathbf v= \Tilde P(i_1,\dots, i_{m-2},\cdot)$
and setting \begin{equation} P(i_1,\dots,i_{m-2},\cdot, \cdot)=
M_{e_1}^T \diag(\mathbf v) M_{e_2}. \label{eq:inddef}\end{equation}

One can check that this definition of $\phi_r$ agrees with our
earlier one, and so is independent of the choice of cherries
defining the sequence $T_2^r,T_3^r,\dots,T_n^r$.

This approach to an inductive definition of $\phi_r$ can be extended
to the case of non-binary trees as follows. For an arbitrary tree
$T^r$, let $\widetilde T^r$ denote any binary tree which resolves
$T^r$, in the sense that $T^r$ can be obtained from $\widetilde T^r$
by contraction of some edges. Extend a choice of parameters $s$ on
$T^r$ to parameters  $\tilde s$ on $\widetilde T^r$ by assigning the
identity matrix to those edges of $\widetilde T^r$ which are
collapsed in $T^r$. Then since $\phi_r(s)=\tilde \phi_r(\tilde s)$,
the inductive definition for binary trees can be applied for
$\widetilde T^r$.

\begin{lem}\label{lem:indphi}
For any $n$-taxon tree $T$, the inductive definition of $\phi_r$
based on Eq. (\ref{eq:inddef}) and outlined above agrees with the
definition in Eq. (\ref{eq:Pdef}).
\end{lem}

\medskip

We also denote by $\phi_r$ the unique extension of this map to a
polynomial map $\phi_r:\C^N\to\C^L$.  The \emph{affine phylogenetic
variety} $V(T)$ for the general Markov model on $T$ is defined as
the closure in $\C^L$ of the image of $\phi_r$.  (Note that this
closure may be taken using either the Zariski topology or the
standard topology on $\C^L$, as the two closures will agree for the
image of a polynomial map.) As has been shown elsewhere
\cite{SSH94,AR03}, this definition is independent of the choice of
the root $r$. $V(T)$ is irreducible, as it is the zero set of a
prime ideal, the kernel of the map between polynomial rings
associated to $\phi_r$.

Now one readily sees the image of $\phi_r$ lies on the hyperplane
defined by the trivial phylogenetic invariant $\sum_{\mathbf i\in
[\kappa]^n} P(\mathbf i)-1=0$. It is therefore natural to pass to the
\emph{projective phylogenetic variety} in $\PP^{L-1}$ by taking a
projective closure. We denote this by $V(T)$ also, making clear by
context whether the affine or projective version is meant.

The \emph{phylogenetic ideals} of all polynomials vanishing on the
affine phylogenetic variety or vanishing on the projective
phylogenetic variety are of course closely related. Generators of
the homogeneous ideal $\aideal_T$ of the projective variety,
together with the trivial invariant, generate the ideal of the
affine variety. Conversely, any homogeneous polynomial in the ideal
of the affine variety is in the homogeneous ideal of the projective
variety. Thus identifying phylogenetic invariants for the general
Markov model means identifying those polynomials vanishing on the
projective phylogenetic variety.

\section{Reparameterization}\label{sec:reparam}

For any projective variety $V\subseteq \PP^m$, let $CV\subseteq
\C^{m+1}$ denote the cone over $V$, that is, the union of the lines
represented by points in $V$. Equivalently, $CV$ is the affine variety
defined by the same polynomials as $V$.

A dense subset of the cone $CV(T)$ admits a parameterization that will
be more useful than the parameterization $\phi_r$ above.  This new
parameterization simplifies many arguments, since it allows matrices
with any row sums to be associated to edges, and no longer requires a
root distribution, or even a specification of a root.

\begin{defn}
Consider an $n$-taxon tree $T$ with $|E|$ edges.  Let $U=\C^K$ with
$K=|E|\kappa^2$. Choose any vertex of $T$ as a root, directing all
edges of $T^r$ away from $r$. View $u\in \C^K$ as a $(2n-3)$-tuple
of complex $\kappa\times\kappa$ matrices $M_e$, one for each edge
$e$ of $T^r$.

Then, in the case that $T$ is binary, let $\psi:U \to \C^L$ be given
inductively as follows, using the sequence of trees
$T^r=T_n^r,T_{n-1}^r,\dots, T_2^r$ chosen in the discussion leading
to Lemma \ref{lem:indphi}:

If $n=2$, $\psi(u)=\psi(M_e)=M_e$, so $\psi$ is the identity map.

If $n>2$, let $\Tilde \psi:\Tilde U\to \C^{\kappa^{n-1}}$ be the map
associated to $T_{n-1}^r$. Then for $u\in U$, define $\tilde u\in
\Tilde U$ by omitting from $u$ the matrices associated to the edges
$e_1,e_2$ of $T_n^r$ not in $T_{n-1}^r$.  Then $\psi(u)=P$, where
$P$ is a $n$-dimensional tensor with 2-dimensional slices given by
first letting $\Tilde P=\tilde \psi(\tilde u)$, $\mathbf v= \Tilde
P(i_1,\dots, i_{n-2},\cdot)$ and setting
$$P(i_1,\dots,i_{n-2},\cdot, \cdot)= M_{e_{1}}^T \diag(\mathbf v)
M_{e_2}.$$

For non-binary trees, modify this construction as indicated for
Lemma \ref{lem:indphi}.
\end{defn}

As in Lemma \ref{lem:indphi}, one sees that this map is independent
of the choice of cherries determining the sequence
$T_2^r,T_3^r,\dots, T_n^r$. Although $\psi$ apparently depends on
the choice of $r$, one can further check that if $r$ is moved from
one vertex of an edge $e$ to the other vertex, we need only
transpose the matrix $M_e$ associated to that edge and the map is
unchanged. Thus the map is independent of the choice of $r$, though
our conception of how components of $\C^K$ are placed into matrices
does depend on $r$. Indeed, all these observations follow from the
observation that $\psi$ can also be defined by a formula like that
in Eq. (\ref{eq:Pdef}), but with the factor
$\boldsymbol{\pi}_r(b_r)$ omitted.

\begin{prop} The closure of $\psi(U)$ in $\C^L$ is the cone $CV(T)$ over
the phylogenetic variety $V(T)$.  \end{prop}

\begin{pf}
To see $\phi_r(S) \subseteq \psi(U)$, suppose
$s=(\boldsymbol{\pi}_r;\{M_e\})\in S$. Let $e_0$ be the one edge of
$T_2^r$, and define $M_{e_0}'=\diag(\boldsymbol{\pi}_r)M_{e_0}$.
With $u=(M_{e_0}', \{M_e\}_{e\ne e_0})$, we find that
$\phi_r(s)=\psi(u)$. Thus $V(T)\subseteq \overline{\psi(U)}$.
Furthermore $\psi(U)$ is a cone, since if $u=(\{M_e\})\in U$ and
$\lambda\in\C$, by picking any particular edge $e_0$ of $T$ and
defining $u'\in U$ to be identical to $u$ but with $\lambda M_{e_0}$
replacing $M_{e_0}$, then $\psi(u')=\lambda\psi(u)$. Thus
$CV(T)\subseteq \overline{\psi(U)}$.

We next show there is a non-empty open, and therefore dense, subset
of $U$ whose image under $\psi$ lies in the cone over $\phi_r(S)$,
and hence in $CV(T)$.  This will imply $ \overline{\psi(U)}
\subseteq CV(T)$.

For simplicity of exposition, assume $T^r$ is binary.

First, if $n=2$, then $\phi_r(S)$ certainly contains those
2-dimensional arrays whose entries add to 1 and none of whose row
sums are 0.  Now the subset of $U$ on which all row sums of
$M_e(=u)$ are non-zero and the total sum of the entries of $M_e(=u)$
is non-zero is an open set. The points in the image under $\psi$ of
this open set lie in the cone over $\phi_r(S)$.

Proceeding inductively, let $e_1$,$e_2$ be the edges of $T_m^r$
which are not in $T_{m-1}^r$, and $e_3$ the third edge meeting them.
We may also suppose $r$ does not lie at the common vertex of
$e_1,e_2,e_3$. Now there is an open $\O_1\subset U$ such that for
points $u\in \O_1$, $M_{e_1}$ and $M_{e_2}$ have all row sums
non-zero.  Letting $D_i$ be the invertible diagonal matrix
constructed from the row sums of $M_{e_i}$, we may write
$$M_{e_i}=D_i M'_{e_i},\ \ i=1,2,$$ where $M'_{e_i}$ has rows summing
to 1.  Let $M_{e_3}'=M_{e_3}D_1D_2$.  Then for any $u\in \O_1$, we
define a new $u'\in\O_1$ as
$$u'=(\{M_e\}_{e\notin \{ e_1,e_2,e_3\}}, M'_{e_1}, M'_{e_2},
M'_{e_3}),$$ so that $\psi(u')=\psi(u)$.  Note that $\omega: \O_1\to
\O_1$ mapping $u\mapsto u'$ is given by rational functions.

Let $\Tilde \psi:\Tilde U\to \C^{\kappa^{m-1}}$ and
$\Tilde\phi_r:\Tilde S\to \C^{\kappa^{m-1}}$ be the
parameterizations associated to $T^r_{m-1}$. Then by induction there
is a non-empty open $\Tilde \O\subset \Tilde U$ such that the image
of all points in $\Tilde \O$ under $\tilde \psi$ lie in the cone
over $\tilde \phi_r(\tilde S)$. Then $\O=\omega^{-1}(\Tilde \O\times
\C^{2\kappa^2} )$ is a non-empty open subset of $U$, and the image
of any point of $\O$ under $\psi$ lies in the cone over $\phi_r(S)$.

If $T^r$ is not binary, slight modifications can be made to the
above argument to obtain the result. \qed\end{pf}

While the definition of $\psi$ has introduced many unnecessary
parameters, in the sense that the dimension of the image is much smaller
than the dimension of the parameter space, it offers us the advantage of
dropping inconvenient requirements --- that row sums of vectors and
matrices be 1 --- that arose from the original probabilistic setting of
the general Markov model.

\section{Flattenings and phylogenetic invariants}\label{sec:SPconj}

To describe the set of phylogenetic invariants we are concerned with, we
require the notion of \emph{flattening} a tensor $P\in \mathbb
\C^{\kappa^n}$ according to an $n$-taxon tree $T$.

Let $e$ be an edge of $T$. Then $e$ induces a split of the taxa
according to the connected components of $T\setminus\{e\}$.  By
reordering the indices in $P$ if necessary, we may assume the split
is $\left \{ \{a_1,\dots, a_k \},\ \{ a_{k+1},\dots, a_n
\}\right\}$. The \emph{flattening of $P$ on $e$} is the
$\kappa^k\times\kappa^{n-k}$ matrix $F=Flat_e(P)$ defined as
follows: Fix any ordering of $J_1=[\kappa]^k$ and
$J_2=[\kappa]^{n-k}$, and for $u\in J_1$, $v\in J_2$, let
$F(u,v)=P(u_1,\dots, u_k,v_1,\dots, v_{n-k})$.

\begin{figure}[h]
\begin{center}
\includegraphics[height=1.25in,width=5in]{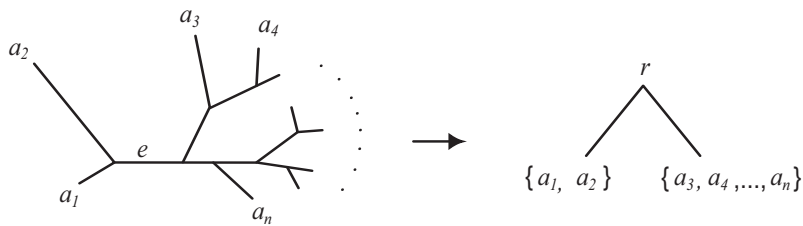}
\end{center}
\caption{Flattening on an edge $e$}\label{fig:edgeflat}
\end{figure}

If the tensor $P=\phi_r(s)$ gives the joint distribution of states
for some parameter choice for the general Markov model on $T$, then
$Flat_e(P)$ can be thought of as a joint distribution for a related
graphical model with less complicated structure: For a tree with at
least 3 leaves, choose the root $r$ to be at one vertex of the edge
$e$, and  imagine at $r$ a $\kappa$-state hidden variable. The
possible joint states at the taxa $a_1,\dots,a_k$ are viewed as a
single $\kappa^k$-state observed variable.  Similarly, the joint
states at the taxa $a_{k+1},\dots, a_n$ are described through a
single $\kappa^{n-k}$-state variable. We thus have a ``coarser''
graphical model with one hidden $\kappa$-state internal node and two
descendent nodes with $\kappa^k$ and $\kappa^{n-k}$ states,
respectively, as depicted in Figure \ref{fig:edgeflat}.  The
flattening of $P$ simply prevents one from examining the finer
structure in the joint distribution array that arises from the
branching of $T$ on either side of $e$.

From this interpretation one readily sees that for any $P\in \phi_r(S)$,
$Flat_e(P)$ has rank at most $\kappa$. Indeed, for the coarser graphical
model, the joint distribution matrix must have the form $$Flat_e(P)=M_1^T
\diag(\pi_r)M_2$$ where $M_1$ and $M_2$ are $\kappa\times\kappa^k$ and
$\kappa\times\kappa^{n-k}$ Markov matrices.

As a result, all $(\kappa+1)\times(\kappa+1)$ minors of $Flat_e(P)$
must vanish. As is classically known, such minors generate the full
prime ideal of polynomials vanishing on matrices of rank $\le
\kappa$, and thus generate all invariants associated to the coarser
model. For the original model on $T$, these minors therefore give
phylogenetic invariants, which we call \emph{edge invariants}
associated to the edge $e$.

We denote by $\mathcal F_{edge}(T)$ the set of all $(\kappa+1)\times
(\kappa+1)$ minors of all flattenings of a $\kappa\times\cdots
\times \kappa$ tensor of $\kappa^n$ indeterminates on edges of $T$.
Of course the choice of ordering of rows and columns in the
flattening introduces factors of $\pm 1$, but as our goal is to
determine ideal generators, we may ignore this issue.

\smallskip

In Section \ref{sec:k2id} we will establish the following, which was
conjectured in \cite{PS}.

\begin{thm}\label{thm:PS} For $\kappa=2$ and any number of taxa $n$, the
phylogenetic ideal $\aideal_T$ for the general Markov model on a
binary $n$-taxon tree $T$ is generated by $\mathcal F_{edge}(T)$,
the $3\times 3$ minors of all edge flattenings of a
$2\times\cdots\times 2$ tensor of indeterminates.  \end{thm}

However, for larger $\kappa$ it is not enough to consider only
2-dimensional edge flattenings (\emph{i.e.}, flattenings to matrices) to
obtain generators of the phylogenetic ideal.  This can be seen already
for the $3$-taxon tree.  In this case, $\mathcal F_{edge}(T)$ is empty,
but for any $\kappa>2$ the phylogenetic ideal contains polynomials of
degree $\kappa+1$ (see \cite{AR03}; for $\kappa=3$ see also \cite{GSS}).
Thus we need at least to consider flattenings of $P$ at internal nodes
of $T$ producing 3-dimensional tensors.

More specifically, consider a trivalent internal vertex $v$ of a
tree $T$, contained in edges $e_1,e_2,e_3$.  Then $v$ induces a
tripartition of the taxa according to the connected components of
$T\setminus\{v,e_1,e_2,e_3\}$. By reordering the indices in $P$ if
necessary, we may assume the tripartition is
$$\left \{ \{a_1,\dots, a_k \},\ \{ a_{k+1},\dots, a_{k+l} \},\
\{a_{k+l+1},\dots, a_n\} \right\}.$$ Then a \emph{flattening of $P$
at $v$} is a $\kappa^k\times\kappa^l\times \kappa^{n-k-l}$ array
$F=Flat_v(P)$ defined as follows: Fix an ordering of
$J_1=[\kappa]^k$, $J_2=[\kappa]^l$, and $J_3=[\kappa]^{n-k-l}$ , and
for $u\in J_1$, $v\in J_2$, $w\in J_3$, let $F(u,v,w)=P(u_1,\dots,
u_k,v_1,\dots, v_{l},w_1,\dots, w_{n-k-l})$.

\begin{figure}[h]
\begin{center}
\includegraphics[height=1.5in, width=5in]{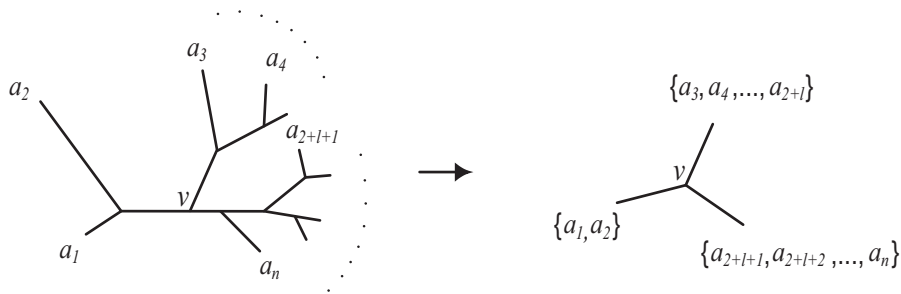}
\end{center}
\caption{Flattening at a vertex $v$}\label{fig:vertexflat}
\end{figure}

As illustrated in Figure \ref{fig:vertexflat}, we think of this
flattening as producing a joint distribution array associated to a
graphical model with one hidden $\kappa$-state internal node and
three descendent nodes with $\kappa^k$, $\kappa^{l}$, and
$\kappa^{n-k-l}$ states, respectively.  Similar to flattenings on
edges, a flattening at an internal node ignores the finer structure
in the joint distribution array that arises from the branching of
$T$ in the three directions leading away from $v$.

An ideal is associated to such a graphical model (1 hidden
$\kappa$-state ancestral node, 3 descendent nodes), and so to the
flattening at a vertex.  While we will investigate such ideals further
in Section \ref{sec:startrees}, already we can formulate a natural
extension of the conjecture of \cite{PS}.

\begin{conj}\label{con:phylotree} For any $\kappa$ and any number of
taxa $n$, the phylogenetic ideal $\aideal_T$ for the general Markov
model on a binary $n$-taxon tree $T$ is the sum of the ideals
associated to the flattenings of $P$ at vertices of $T$.  \end{conj}

That this conjecture is identical to Theorem \ref{thm:PS} when
$\kappa=2$ follows from work of Landsberg and Manivel \cite{lands}.
They show that in this special case the ideal associated to a vertex
flattening is the sum of those associated to the edge flattenings on
the three edges containing the vertex.  (The Landsberg-Manivel
result is a special case of a conjecture in \cite{GSS}. We will give
a new and simpler proof of this case, and several additional cases,
as Corollary \ref{cor:GSS345}.)

Of course the notion of flattening at a vertex can be extended in a
straightforward way for vertices of valence $>3$, and the conjecture
formulated for non-binary trees as well. The extended conjecture for
non-binary trees remains open even for $\kappa=2$.

\medskip

Although we will primarily need to refer to the 2- and 3-dimensional
flattenings of a tensor $P$ on an edge or at a vertex of a tree
$T$, the notion naturally extends to flattenings based on any partition
of the set of labels (taxa) associated to the indices of $P$. For
instance, an $n$-dimensional $\kappa\times\cdots\times \kappa$ tensor
$P$ with associated labels $a_1,\dots, a_n$ can be flattened according
to the partition $\left \{ \{a_1\},\dots, \{a_{n-2}\},\{a_{n-1},a_n\}
\right \}$ to give an $(n-1)$-dimensional $\kappa\times\cdots
\times\kappa\times\kappa^2$ tensor.  We use such a flattening, where
$a_{n-1},a_n$ are in a cherry, in Section \ref{sec:k2id}. Flattenings
according to arbitrary bipartitions also appear in Section
\ref{sec:startrees}.

\section{The algebra of tensors, trees, and parameters}\label{sec:alg}

In this section we define binary operations on trees, model parameters
on trees, and tensors. These operations, all denoted by the same symbol
`$\st$', exhibit relationships that will make them useful in later sections.

\smallskip

\noindent {\bf Tensors:} If $Q$ and $R$ are $m$- and $n$-dimensional
tensors of `matching size $\kappa$' in the last and first index
respectively, then we define an $l=(m+n-2)$-dimensional tensor $Q\st R$
by $$(Q\st R) (i_1,\dots i_l) =\sum_{j=1}^\kappa
Q(i_1,\dots,i_{m-1},j)R(j,i_{m},\dots,i_{l}).$$ For $m=n=2$, this is of
course just matrix multiplication.

More generally, if the $p$th index of $Q$ and the $q$th index of $R$
both run through $[\kappa]$, we may define $Q\st_{p,q} R$ by a similar
sum. However, to keep our notation less cumbersome, we will generally
try to express products using the last and first indices.

\smallskip

\begin{figure}[h]
\begin{center}
\includegraphics[height=1.2in, width=5in]{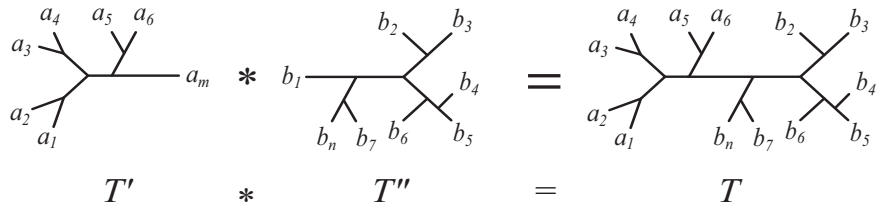}
\end{center}
\caption{The $\st$ operation on trees}\label{fig:treestar}
\end{figure}

\noindent {\bf Trees:} Suppose $T'$ is a tree with taxa $a_1,a_2,\dots,
a_m$, and $T''$ is a tree with taxa $b_1,b_2,\dots, b_n$.  Then by
$T'\st T''$ we mean the $(m+n-2)$-taxon tree with taxa $a_1,\dots,
a_{m-1},b_2,\dots, b_n$ obtained by first identifying the vertices $a_m$
and $b_1$, and then deleting this vertex, replacing the two edges it
lies in with a single conjoined edge, as illustrated in Figure
\ref{fig:treestar}.

\noindent {\bf Parameters:} Consider trees $T'$, $T''$, and $T=T'\st
T''$ with $m$, $n$, and $m+n-2$ taxa.  Then from Section
\ref{sec:reparam} we have the parameterizations \begin{align*}
\psi'&:U'\to \C^{\kappa^m},\\ \psi''&:U''\to \C^{\kappa^n},\\ \psi &:
U\to \C^{\kappa^{m+n-2}}, \end{align*} of the cones over the associated
phylogenetic varieties.

To impose directions on the edges of the trees for notational purposes,
root $T'$ and $T$ at $a_1$, and $T''$ at $b_1$.  Then for $u'\in U'$,
$u''\in U''$, we define $u'\st u''\in U$ by retaining for each edge of
$T$ except the conjoined one the matrix associated to the edge in either
$u'$ or $u''$, and for the conjoined edge using the product of the
matrices in $u'$ and $u''$ associated to its parts.

\smallskip

One readily sees that these three definitions imply the following.

\begin{lem}\label{lem:starparam} $\psi(u'\st u'')=\psi'(u')\st
\psi''(u'')$. \end{lem}

\begin{lem}\label{lem:closst} If $T=T'\st T''$, then
$CV(T)=\overline{CV(T')\st CV(T'')}.$ \end{lem} \begin{pf} It is
clear that $$U=U'\st U''=\{u'\st u'' ~|~ u'\in U',\ u''\in U''\}.$$
Thus by Lemma \ref{lem:starparam}, $$CV(T)=\overline{\psi(U)} =
\overline{\psi'(U') \st \psi''(U'')} = \overline{CV(T')\st
CV(T'')}.$$ \qed\end{pf}

This result will be strengthened in Corollary \ref{cor:pvalg}.

\smallskip

In the special case when $T''$ is a 2-taxon tree, $T'\st T''$ is
isomorphic to $T'$.  Then $u''=\psi(u'')$ is simply a
$\kappa\times\kappa$ matrix.  Informally, one can think of $\psi'(u')\st
\psi''(u'')$ as the result of `extending' the edge of $T'$ terminating
at $a_m$ and associating to the edge extension the matrix $u''$.

Considering invertible matrices $u''$, we get an action of
$GL(\kappa,\C)$ on both $U'$ and $\psi'(U')$. Thus $GL(\kappa,\C)$ acts
on the closure, $CV(T')$, as well.  Viewing the action described here as
operating in `the last index' of a tensor in $V_{T'}$, we similarly have
an action in the other indices.  These actions of
$GL(\kappa,\C)$ are of course just restrictions of the natural actions
of that group on the set of all $\kappa\times\cdots \times\kappa$
tensors: For $j=1,\dots,n$, the `$j$th index' action is defined by
$P\mapsto P\st_{j,1}A$ for $A\in GL(\kappa,\C)$.

\section{Models on star trees}\label{sec:startrees}

In this section, we step back from the phylogenetic tree setting,
and consider in more depth the hidden naive Bayes models of
\cite{GSS}. Most of our results will be needed for application to
phylogenetic varieties.  However, we develop this material in
slightly greater generality than we need for phylogenetic
applications, and so obtain partial results on a conjecture of
\cite{GSS} as well.

\begin{figure}[h]
\begin{center}
\includegraphics[height=1.2in, width=5in]{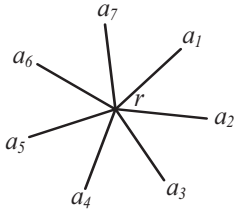}
\end{center}
\caption{Graphical depiction of a hidden naive Bayes model}\label{fig:hnbm}
\end{figure}

The graphical models of this section are based on a star tree, as in
Figure \ref{fig:hnbm}, with one internal vertex $r$, connected by
edges to $n$ leaves $a_1,a_2,\dots, a_n$.  A hidden random variable
associated to $r$ has $\kappa$ possible states, with probability
distribution given by a vector $\boldsymbol{\pi}_r$.  Each leaf
$a_i$ has associated to it a random variable with $l_i$ states, and
Markov matrices $M_i$ of size $\kappa\times l_i$ give conditional
probabilities of observing the various states at $a_i$ given the
state at $r$.

As in the phylogenetic situation, such a model defines a projective
variety, the closure of the set of joint distributions of
observations at the leaves arising from this parameterization.  We
denote this variety by $V(\kappa;l_1,l_2,\dots, l_n)$, and the
homogeneous ideal defining it by $\aideal(\kappa;l_1,l_2,\dots,
l_n)$.

As pointed out in \cite{GSS}, the variety
$V(\kappa;l_1,l_2,\dots,l_n)$ can be viewed more geometrically as
the $\kappa$-secant variety of the Segre product $\PP^{l_1-1}\times
\PP^{l_2-1}\times\dots\times\PP^{l_n-1}$. Here `$\kappa$-secant
variety' means the closure of the union of the
$(\kappa-1)$-dimensional affine spaces spanned by collections of
$\kappa$ points on the original variety, so, for instance, the
2-secant variety arises from points on secant lines.

Note that $V(\kappa;\kappa,\kappa,\kappa)= V(T_3)$, the phylogenetic
variety for a $\kappa$-state, 3-taxon tree. The varieties
$V(\kappa;\kappa^k,\kappa^l,\kappa^{n-k-l})$, with $k,l,n-k-l\ge 1$, are
the ones that arose in Section \ref{sec:SPconj}, in the discussion of
flattenings of tensors at vertices of phylogenetic trees. Moreover,
flattenings on edges involve $V(\kappa;\kappa^k,\kappa^{n-k}$), the
variety of rank $\kappa$ matrices of size $\kappa^k\times\kappa^{n-k}$,
which is well understood classically.

\smallskip

Our first goals are to show Theorems \ref{thm:ststar} and
\ref{thm:idstar}: Given a set $\mathcal F$ of polynomials
set-theoretically (respectively, ideal-theoretically) defining
$V(\kappa;\kappa,\kappa,\dots,\kappa)$ for the $n$-leaf star tree,
then for any $l_i\ge \kappa$ we can explicitly construct polynomials
set-theoretically (respectively, ideal-theoretically) defining
$V(\kappa;l_1,l_2,\dots,l_n)$.

Previous to these theorems, we know of only one general result
concerning defining polynomials of $V(\kappa;l_1,l_2,\dots,l_n)$:
When $\kappa=2$, for any number of leaves, \cite{lands} gives a
natural set of polynomials defining the variety as a set.

\smallskip

For our application to phylogenetic trees, the assumption that internal
nodes are trivalent means only the case $n=3$ is needed.  We therefore
summarize known results on $V(\kappa;\kappa,\kappa,\kappa)=V(T_3)$ for
small $\kappa$.

For $\kappa=2$, as noted in \cite{GSS,AR03,lands}, $V(T_3)=\PP^7$,
and so $\{0\}$ is the full prime ideal defining the variety.

For $\kappa=3$, a generating set $\mathcal F$ for the prime ideal
may be taken to be the 27 quartic polynomials in \cite{GSS}, first
found in \cite{MR85b:15039} but also obtained from the construction
in \cite{AR03}.

For $\kappa\ge 4$, finding an explicit set $\mathcal F$ that even
set-theoretically defines the variety is still an open problem.
However, any polynomial vanishing on the variety must be of degree
at least $\kappa+1$.

When $\kappa=4$, all degree 5 polynomials vanishing on the variety
form an explicitly-known 1728-dimensional vector space.  This
dimension is computed in \cite{HagCRM,lands}, and an explicit
construction for general $\kappa$ is given in \cite{AR03} that
produces a spanning set when $\kappa=4$.  Moreover, off another
explicitly-known variety, the vanishing of these polynomials does
distinguish points of $V(T_3)$. However, an explicit degree 9
polynomial is known which vanishes on $V(4;3,3,3)$ (see \cite{GSS}
for a statement, or \cite{MR85b:15039} for the construction), and
from this polynomial one can obtain degree 9 polynomials vanishing
on $V(4;4,4,4)=V(T_3)$ by evaluation on $3\times 3\times 3$
subarrays of a $4\times 4\times 4$ tensor. By consideration of the
multidegrees of each monomial term in its many variables, one can
show that these degree 9 polynomials cannot be generated by the
degree 5 ones.

We also note that if $\phi$ is the parameterization arising from the
general Markov model on $T_3$, then for all $\kappa\ge 2$ the image
of $\phi$ is strictly smaller than its closure $V(T_3)$. This is
pointed out in Section 9 of \cite{AR03}, but in the terminology of
\cite{MR85b:15039} is simply the statement that `rank $\kappa$' is a
strictly stronger statement than `border rank $\kappa$' for
$\kappa\times\kappa\times\kappa$ tensors.

\medskip

By modifying the approach of Section $\ref{sec:reparam}$, it is possible
to parameterize a dense subset of the cone $CV(\kappa;l_1,l_2,\dots,
l_n)$ using parameters which are arbitrary matrices. We leave the
details to the reader, but denote this parameterization by
$\psi_{\kappa;l_1,\dots, l_n}$, where $$\psi_{\kappa;l_1,\dots,
l_n}:U_{\kappa;l_1,\dots, l_n}\to \C^L,\ \ U_{\kappa;l_1,\dots,
l_n}=\C^{\kappa(l_1+\cdots+l_n)},\ L=l_1 l_2\cdots l_n,$$ and if
$P=\psi(M_1,M_2,\dots, M_n)$ then $$P(i_1,\dots, i_n)=\sum_{k=1}^\kappa
\prod_{j=1}^n M_j(k,i_j).$$ Here $M_j\in M(\kappa,l_j,\C)$, the set of
complex $\kappa\times l_j$ matrices.

\smallskip
In order to relate $V(\kappa;l_1,l_2,\dots, l_n)$ to
$V(\kappa;\kappa,\kappa,\dots, \kappa)$ we need the following lemma.
It can be interpreted as describing the effect of extending one edge
of the star tree, and associating a (non-square) matrix to that
extension, as was explained at the end of Section \ref{sec:alg}.

\begin{lem} \label{lem:extedge} Let $P\in CV(\kappa;l_1,l_2,\dots,l_n)$
and let $A\in M (l_n, l'_n,\C)$. Then $A$ defines a map
$CV(\kappa;l_1,l_2,\dots, l_n)\to CV(\kappa;l_1,l_2,\dots, l'_n)$ by
$P\mapsto P*A$.  Furthermore, \begin{enumerate} \item[\it (i)] If $rank(A)=l_n'$,
then $CV(\kappa;l_1,l_2,\dots, l_n)*A$ is dense in\\
$CV(\kappa;l_1,l_2,\dots, l'_n)$.  \item[\it (ii)] If $\kappa\le l_n$ then
$CV(\kappa;l_1,l_2,\dots, l_n)*M(l_n,l'_n,\C)$ is dense in\\
$CV(\kappa;l_1,l_2,\dots, l'_n)$.  \end{enumerate} \end{lem}

\begin{pf} Suppose first that
$P=\psi_{\kappa;l_1,\dots,l_n}(M_1,M_2,\dots, M_n)$, with complex
$\kappa\times l_i$ matrix parameters $M_i$, $i=1,2,\dots,n$ associated
to the $n$ edges of $T$ directed away from the internal node. Then $$P
\st A = \psi_{ \kappa ; l_1, \dots, l_{n-1},l_n'} (M_1,M_2,\dots,
M_{n-1}, M_nA),$$ hence $P\st A \in CV(\kappa;l_1,l_2,\dots,l'_n)$.
Since $P\st A \in CV(\kappa;l_1,l_2,\dots,l'_n)$ for $P$ in a dense
subset of $CV(\kappa;l_1,l_2,\dots,l_n)$, it follows that $P\st A \in
CV(\kappa;l_1,l_2,\dots,l'_n)$ for all $P\in
CV(\kappa;l_1,l_2,\dots,l_n)$.

\smallskip

Now suppose $\rank(A)=l'_n$. Then as $M_n$ ranges through all
$\kappa\times l_n$ complex matrices, $M_nA$ ranges through all
$\kappa\times l_n'$ complex matrices. Thus $$\psi_{\kappa;l_1,\dots,
l_n}(U_{\kappa;l_1,\dots, l_n})\st
A=\psi_{\kappa;l_1,\dots,l_{n-1},l'_n}(U_{\kappa;l_1,\dots, l'_n}),$$
and so a subset of $CV(\kappa;l_1,\dots,l_n)\st A$ is dense in
$CV(\kappa;l_1,\dots,l'_n)$.

Finally suppose $\kappa\le l_n$. Then as $M_n$ ranges through all
$\kappa\times l_n$ complex matrices and $A$ through all $l_n\times l_n'$
matrices, $M_nA$ ranges through all $\kappa\times l_n'$ complex
matrices. Thus $$\psi_{\kappa;l_1,\dots, l_n}(U_{\kappa;l_1,\dots, l_n})
\st M(l_n,l'_n,\C) = \psi_{\kappa;l_1,\dots,l_{n-1},l'_n}
(U_{\kappa;l_1,\dots, l'_n}).$$ Therefore a subset of
$CV(\kappa;l_1,\dots,l_n)\st M(l_n,l'_n,\C)$ is dense in\\
$CV(\kappa;l_1,\dots,l'_n)$.  \qed\end{pf}

\begin{rem} For non-zero $P$ and $A$ as in the proof, it is possible
for $P\st A$ to be a zero tensor. Thus while the above lemma could
be formulated in terms of a rational map between the underlying
projective varieties, it is slightly easier for us to consider a
polynomial map on the cones.  \end{rem}

By permuting indices Lemma \ref{lem:extedge} can be applied in any
index, not just the last. As shorthand, we will refer to this as
letting an $l_k\times l_k'$ matrix `act in the $k$th index.' By
considering only invertible $l_k\times l_k$ matrices, we have a
group action of $GL(l_k,\C)$ in the $k$th index, and so an action of
$GL(l_1,\C)\times \cdots\times GL(l_n,\C)$ on
$V(\kappa;l_1,\dots,l_n)$. While this group action underlies the
dimension computations of \cite{lands}, our work will emphasize the
utility of non-square and non-invertible matrices as well.

\begin{thm} \label{thm:ststar} Consider an $n$-leaf star tree.
Suppose $l_1,l_2,\dots,l_n \ge \kappa$. Let $\mathcal F$ be any set
of polynomials whose zero set is
$V(\kappa;\kappa,\kappa,\dots,\kappa)$. For $k=1,2,\dots, n$, let
$Z_k=(z_{ij}^k)$ be $l_k\times\kappa$ matrices of indeterminates.
For an $l_1\times l_2\times\dots\times l_n$ tensor $P$ of
indeterminates, let $\Tilde P$ be the
$\kappa\times\kappa\times\dots\times \kappa$ tensor that results
from letting each $Z_k$ act formally in the $k$th index of $P$. Let
$\Tilde{ \mathcal F}$ denote the set of polynomials in the entries
of $P$ obtained from those in $\mathcal F$ by substituting into them
the entries of $\Tilde P$, expressing the results as polynomials in
the $z_{ij}^k$, and then extracting the coefficients. Let $\mathcal
F_{edge}$ denote the set of $(\kappa+1)\times(\kappa+1)$ minors of
the $n$ flattenings of $P$ on edges of the star tree. Finally, let
$\mathcal F(\kappa;l_1,l_2,\dots, l_n)= \Tilde { \mathcal F} \cup
\mathcal F_{edge}$.

Then $\mathcal F(\kappa;l_1,l_2,\dots, l_n)$ defines
$V(\kappa;l_1,l_2,\dots, l_n)$ set-theoretically. \end{thm}

\begin{pf} We first observe that all polynomials in $\mathcal
F(\kappa;l_1,l_2,\dots,l_n)$ vanish on the cone
$CV(\kappa;l_1,l_2,\dots,l_n)$: Polynomials in $\mathcal F_{edge}$
must vanish there, since the model has $\kappa$ states at the
internal node, so all 2-dimensional flattenings on edges must have
rank $\le \kappa$ on the parameterized subset of the variety, and
hence on the whole variety. Polynomials in $\Tilde {\mathcal F}$
must vanish there, since for all assignments of values to the
$z_{ij}^k$, if $P\in CV(\kappa,l_1,\dots, l_n)$ then, by Lemma
\ref{lem:extedge}, $\Tilde P\in CV(\kappa;\kappa,\dots ,\kappa )$.

Now suppose all polynomials in $\mathcal
F(\kappa;l_1,l_2,\dots,l_n)$ vanish on a tensor $P_0\in
\C^{l_1l_2\cdots l_n}$.  Then, flattening $P_0$ on the edge of the
tree leading to $a_n$ gives a matrix of rank $l\le \kappa$, so we
can write $$P_0=Q'_0*B_n'$$ where $Q_0'$ is an $l_1\times l_2\times
\dots \times l_{n-1}\times l$ tensor and $B_n'$ is an $l \times l_n$
matrix. Construct a $\kappa\times l_n$ matrix $B_n$ of rank $\kappa$
by augmenting $B_n'$ with additional rows. Similarly augment $Q_0'$
with additional zero entries to obtain an $l_1\times l_2\times \dots
\times l_{n-1}\times \kappa$ tensor $Q_0$ with $P_0=Q_0*B_n$.  Now
there exists an $l_n\times \kappa$ matrix $A_n$ so that $B_nA_n=I$,
the identity matrix.  Thus $P_0\st A_n\st
B_n=Q_0*B_n*A_n*B_n=Q_0*I*B_n=P_0$.

Proceeding similarly for the other taxa, we obtain matrices $A_k$,
$B_k$ such that
$$(P_0*_{k,1}A_k)\st_{k,1}B_k=P_0*_{k,1}(A_k*B_k)=P_0.$$ By
simultaneously letting each $A_k$ act in the $k$th index of $P_0$,
we obtain a $\kappa\times\kappa\times\dots\times\kappa$ tensor
$\tilde P_0$.  Because all polynomials in $\Tilde {\mathcal F}$
vanish on $P_0$, all polynomials in $\mathcal F$ vanish on $\tilde
P_0$. Thus by our choice of $\mathcal F$, $\tilde P_0\in CV(\kappa;
\kappa, \kappa, \dots, \kappa)$.  Since, by repeated applications of
Lemma \ref{lem:extedge}, letting each $B_k$ act in the $k$th index
maps $CV(\kappa; \kappa, \kappa, \dots, \kappa)$ to $CV(\kappa; l_1,
l_2, \dots,l_n)$, and maps $\tilde P_0$ to $P_0$, we see $P_0\in
CV(\kappa; l_1,l_2,\dots,l_n)$. \qed\end{pf}

We now state an ideal-theoretic version of this result.

\begin{thm}\label{thm:idstar} Suppose $l_1,l_2,\dots,l_n \ge \kappa$,
and $\mathcal F$ is a set of polynomials generating
$\aideal(\kappa;\kappa,\kappa,\dots,\kappa)$.  Then the set
$\mathcal F(\kappa;l_1,l_2,\dots, l_n)$ constructed from $\mathcal
F$ as in Theorem \ref{thm:ststar} generates
$\aideal(\kappa;l_1,l_2,\dots,l_n)$.
\end{thm}

Since the key argument in the proof of Theorem \ref{thm:idstar} will
be used again in Section \ref{sec:k2id}, we present it as a lemma.

\begin{lem} \label{lem:reparg} Let $V_1$ and $V_2$ be subvarieties of
$\C^{n m_1}$ and $\C^{nm_2}$, respectively, with $m_1\le m_2$, such that,
when points are written as $n\times m_1$ and $n\times m_2$ matrices,
$$V_1=V_1\st M(m_1,m_1,\C),$$ and $$V_2=\overline{V_1\st
M(m_1,m_2,\C)}.$$ Let $\aideal_i$ denote the ideal of all
polynomials vanishing on $V_i$.

Then $\aideal_2$ is generated by the $(m_1+1) \times (m_1+1)$ minors
of an $n\times m_2$ matrix $P$ of indeterminates, together with all
polynomials of the form $f(P\st A)$, where $f\in\aideal_1$ and $A\in
M(m_2,m_1,\C)$.  \end{lem}

\begin{pf} Let $\bid$ denote the ideal generated by the
$(m_1+1)\times (m_1+1)$ minors, together with the polynomials $f(P\st
A)$ described above.

First we show $\aideal_2\supseteq \bid$.  It is enough to show the
specified generators of $\bid$ vanish on $V_1\st M(m_1,m_2,\C)$.
Since all points in this set are matrices of rank at most $m_1$, the
specified minors vanish there.  To see the $f(P\st A)$ vanish there,
consider a point $Q_0\st B$ where $Q_0\in V_1$, $B\in
M(m_1,m_2,\C)$.  Then $Q_0\st B \st A \in V_1$ since $B\st A \in
M(m_1,m_1,\C)$. Thus $f(P\st A)$ vanishes at $Q_0\st B$.

\smallskip

Our argument that $\aideal_2\subseteq\bid$ is more involved.

Note $GL(m_2,\C)$ acts on $V_1\st M(m_1,m_2,\C)$, and hence on $V_2$
as well.  Consider the degree $m$ homogeneous component
$\aideal_2^{(m)}$ of $\aideal_2$. Then the $GL(m_2,\C)$-action on
$V_2$ gives a representation of $GL(m_2,\C)$ on $\aideal_2^{(m)}$,
in which $C\in GL(m_2,\C)$ maps the polynomial $g(P)\mapsto g(P\st
C)$. Since $GL(m_2,\C)$ is reductive, this representation decomposes
into a sum of irreducible ones. Consider now one of the irreducible
subspaces, $W$. It will be enough to show that $W\subseteq\bid$.

Consider any non-zero polynomial $g(P)\in W$. Let $Q$ denote a
$n\times m_1$ matrix of indeterminates. Then for any $B\in
M(m_1,m_2,\C)$, the polynomial $g_B(Q)=g(Q\st B)$ vanishes on $V_1$,
since $Q_0\mapsto Q_0\st B$ maps $V_1$ to $V_2$. Thus $g_B\in
\aideal_1$.

Suppose first that for all $B\in M(m_1,m_2,\C)$ the polynomial $g_B(Q)$
is identically zero.  Then $g$ must vanish on all $n\times m_2$ matrices
of rank at most $m_1$, since any such matrix can be written as $Q_0\st
B$ for some complex matrices $Q_0\in M(n,m_1,\C)$, $B\in M(m_1,m_2,\C)$,
and then $g(Q_0\st B)=g_B(Q_0)=0$.  Thus if all $g_B$ are identically
zero, then $g$ is in the ideal generated by $(m_1+1) \times(m_1+1)$
minors of $P$, and hence $g \in \bid$.

Suppose, then, that for some $B$ the polynomial $g_B$ is not identically
zero. Let $D\in M(m_2,m_1,\C)$ be chosen so that $h(P)=g_B(P\st D)$ is a
non-zero polynomial. Such a $D$ must exist since $m_1\le m_2$.  (For
instance, $D$ may be taken so that its first $m_1$ rows form an identity
and the remaining rows are zero.) Then $h(P)=g(P\st DB)$, where $DB$ is
a complex $m_2\times m_2$ matrix that is generally not invertible.

Nonetheless, the irreducibility of $W$ implies that $h(P)\in W$.
This is simply because $W$ is closed in $\aideal_2^{(m)}$, and so
must contain the closure of the orbit of $g$ under $GL(m_2,\C)$, and
this closure contains $g(P\st DB)$.

Now since $g(P)\in W$, $h(P)=g_B(P\st D)\in W$, and $g_B\in
\aideal_1$, the irreducibility of $W$ implies $g(P)$ is in the span
of polynomials of the form $f(P\st A)$ where $f\in \aideal_1$ and
$A\in M(m_2,m_1,\C)$. Thus in this case as well, $g\in\bid$.
\qed\end{pf}

\begin{pf}[Proof of Theorem \ref{thm:idstar}] Let $\aideal = \aideal(
\kappa; l_1, \dots, l_n)$, and let $\bid$ be the ideal generated by
$\mathcal F(\kappa;l_1,\dots, l_n)$, the set defined in Theorem
\ref{thm:ststar}. Note that $\bid$ is equivalently described as
generated by $\mathcal F_{edge}\cup \Tilde{\Tilde{\mathcal F}}$,
where $\Tilde{\Tilde{\mathcal F}}$ denotes the set of all
polynomials of the form $f(\Tilde{P})$ where $f \in \mathcal F$ and
$\Tilde P$ is obtained from a tensor $P$ of indeterminates by the
action of numerical matrices $Z_k\in M(l_k,\kappa,\C)$ in each index
$k$.

That $\aideal\supseteq\bid$ was shown in the proof of Theorem
\ref{thm:ststar}.  To establish $\aideal\subseteq\bid$.  we proceed
by induction on the number of indices $k$ such that $l_k>\kappa$,
the base case of zero being trivial.

If at least one such $l_k>\kappa$ exists, we may assume
$l_n>\kappa$. Then let $V_1=CV(\kappa; l_1,\dots, l_{n-1},\kappa)$
and $V_2=CV(\kappa; l_1,l_2,\dots, l_n)$.  We view points on $V_1$
and $V_2$ as $l_1\cdots l_{n-1}\times \kappa$ and $l_1\cdots
l_{n-1}\times l_n$ matrices, respectively, by flattening on the edge
of the star tree leading to the $n$th leaf.  Using Lemma
\ref{lem:extedge} we see that $V_1*M(\kappa,\kappa,\C)=V_1$ and,
since $l_n>\kappa$, that $V_2=\overline{V_1*M(\kappa,l_n)}$.
Therefore we may apply Lemma \ref{lem:reparg}, and obtain that
$\aideal$ is generated by the $(\kappa+1)\times(\kappa+1)$ minors of
the flattening of $P$ on the edge to the $n$th leaf, together with
all polynomials $f(P\st A)$ where $f\in\aideal(\kappa; l_1,\dots,
l_{n-1},\kappa)$ and $A\in M(l_n,\kappa,\C)$. We thus need only show
such $f(P\st A)$ are in $\bid$.

Now by induction, $\aideal(\kappa;l_1,\dots, l_{n-1},\kappa)$ is
generated by $(\kappa+1)\times (\kappa+1)$ minors of edge
flattenings of an $l_1\times\cdots\times l_{n-1}\times\kappa$ tensor
$Q$ of indeterminates, together with polynomials of the form
$h(\Tilde Q)$, where $h\in \aideal(\kappa;\kappa,\dots,\kappa)$ and
$\Tilde Q$ is a $\kappa\times\cdots \times\kappa$ tensor obtained
from $Q$ by letting elements of $M(l_i,\kappa,\C)$ (respectively
$M(\kappa,\kappa,\C)$) act on $Q$ in the $i$th index for each $i\ne
n$ (respectively $i=n$). We may thus assume $f$ itself has one of
these forms.

In the first case, where $f\in\aideal(\kappa; l_1,\dots,
l_{n-1},\kappa)$ is a minor of an edge flattening for the model, we
see $f$ vanishes on all tensors $Q$ that have rank at most $\kappa$
when flattened on a certain edge $e$ not leading to the $n$th leaf.
But if $P$ is an $l_1\times\cdots\times l_n$ tensor with
$\rank(Flat_e(P))\le\kappa$, then $\rank(Flat_e(P*A)) \le \kappa$ as
well, for all $A\in M(l_n,\kappa,\C)$.  Thus $f(P*A)$ vanishes on
all tensors such that $\rank(Flat_e(P))\le\kappa$, and so $f(P\st
A)$ is in the ideal generated by $(\kappa+1)\times(\kappa+1)$ minors
from edge flattenings of $P$.

In the second case, where $f=h(\Tilde Q)$, we find $f(P*A)=h(\Tilde P)$
where $\tilde P$ is obtained from $P$ by letting elements of
$M(l_i,\kappa,\C)$ act on $P$ in the $i$th index for each $i$, and $h$
vanishes on $V(\kappa;\kappa,\cdots,\kappa)$.

Thus in either case $f(P\st A)\in \bid$.  \qed\end{pf}

\begin{rem}
It is natural to ask whether a smaller set of polynomials than the
set described above --- namely, those constructed by evaluation of
polynomials in $\mathcal F$ on all $\kappa\times\cdots\times\kappa$
subarrays of a $l_1\times\cdots\times l_n$ array of indeterminates
--- is sufficient to define the variety $V(\kappa;l_1,\dots,l_n)$.
Indeed, Lemma \ref{lem:ppgss} below shows it does in the special
case $\kappa=2$, assuming elements of $\mathcal F$ have a special
form.

However, in general this subset does not even define the variety
set-theoretically. To see this, consider the $3\times 3\times 4$
tensor
$$P=e_1\otimes e_1\otimes f_1+
e_2\otimes e_2\otimes f_2+ e_3\otimes e_3\otimes f_3+ e_1\otimes
e_2\otimes f_4,$$ where the $e_i$  are the standard basis vectors
for $\C^3$ and the $f_i$ the standard basis vectors for $\C^4$. That
all $3\times 3 \times 3$ subarrays of $P$ are in $V(3;3,3,3)$ is
clear from the form of $P$. One can verify that $P\notin V(3;3,3,4)$
by checking the non-vanishing at $P$ of some of the polynomials
constructed in Theorem \ref{thm:idstar}.

\end{rem}

As a corollary to Theorem \ref{thm:idstar}, we prove several cases
of Conjecture 21 in \cite{GSS} on the ideals
$\aideal(2;l_1,\dots,l_n)$. We note the $n=3$ case was first proved
in \cite{lands} by invoking sophisticated methods of Weyman
\cite{MR2004d:13020}.

\begin{cor} \label{cor:GSS345} For $n\le 5$, the ideal
$\aideal(2;l_1,\dots,l_n)$ associated to the hidden naive Bayes
model with a 2-state hidden variable and $n$ observed variables with
$l_1,\dots,l_n$ states, is generated by the $3\times 3$ minors of
all 2-dimensional flattenings associated to bipartitions of the
observed variables.  \end{cor}

\begin{pf} Since there are no polynomials vanishing on
$V(2;2,2,2)=\PP^7$, by Theorem \ref{thm:idstar} the set of polynomials
vanishing on $V(2;l_1,l_2,l_3)$ is generated by edge invariants.

By calculations of \cite{GSS}, the statement holds for the two cases
$V(2;2,2,2,2)$ and $V(2;2,2,2,2,2)$. The corollary then follows from
Lemma \ref{lem:ppgss} below.  \qed\end{pf}

\begin{lem}\label{lem:ppgss} Suppose, for the $n$-leaf star tree, that
the ideal $\aideal(2;2,\dots,2)$ is generated by the $3\times 3$
minors of all 2-dimensional flattenings of $2\times\cdots\times 2$
tensors according to bipartitions of the observed variables. Then
$\aideal(2;l_1,\dots,l_n)$ is generated by the $3\times 3$ minors of
all 2-dimensional flattenings of $l_1\times\cdots\times l_n$ tensors
according to bipartitions of the observed variables.  \end{lem}

\begin{pf} By Theorem \ref{thm:idstar}, $\aideal(2;l_1,\dots,l_n)$ is
generated by all $3\times 3$ minors of edge flattenings of an
$l_1\times\cdots\times l_n$ tensor of indeterminates $P$, together
with all $3\times 3$ minors of all 2-dimensional flattenings of all
$\Tilde P$, where $\Tilde P$ denotes a $2\times\cdots \times 2$
tensor obtained from  $P$ by an action in each index $i$ by matrices
$A_i\in M(l_i,2,\C)$.  One readily sees such flattenings of $\Tilde
P$ can be expressed as $\Tilde F=B_1*F*B_2$, where $F$ is the
corresponding flattening of $P$ and the $B_j$ are matrices depending
on the $A_i$. But then the $3\times 3$ minors of such a flattening
of $\Tilde P$ will be zero provided $F$ has rank $\le 2$. Thus these
polynomials are in the ideal $\bid$ generated by $3\times 3$ minors
of flattenings of $P$.  Therefore $\aideal(2;l_1,\dots,l_n)\subseteq
\bid$.

That $\aideal(2;l_1,\dots,l_n)\supseteq \bid$ is clear.
\qed\end{pf}

A proof of the full Conjecture 21 of \cite{GSS} will therefore
follow from the following special cases:

\begin{conj}\label{con:GSSsimp}(Garcia,Stillman,Sturmfels) The ideal
$\aideal(2;2,2,\dots,2)$, that is, the ideal associated to the
hidden naive Bayes model with a 2-state hidden variable and $n$
2-state observed variables, is generated by the $3\times 3$ minors
of all 2-dimensional flattenings arising from bipartitions of the
observed variables.
\end{conj}

\section{Set-theoretic description of the phylogenetic variety:
arbitrary $\kappa$.}\label{sec:kanyvar}

For the remainder of this paper, we return to the consideration of
models on phylogenetic trees.  We first establish a set-theoretic result
that provides some evidence for Conjecture \ref{con:phylotree}, for
arbitrary $\kappa$.

\begin{thm} \label{thm:kany}  For an $n$-taxon tree $T$, let
$\mathcal F(T)$ be the union of all sets of polynomials $\mathcal
F(\kappa; l_1,l_2,\dots l_n)$, defined as in Theorem
\ref{thm:ststar}, associated to flattenings at nodes of $T$. Then
the zero set of $\mathcal F(T)$ is the phylogenetic variety $V(T)$.
\end{thm}

More informally, in conjunction with Theorem \ref{thm:ststar} this
means that from polynomials whose zero set is
$V(\kappa;\kappa,\dots,\kappa)$ one can explicitly construct
polynomials whose zero set is $V(T)$ for any $n$-taxon tree $T$.

In particular, knowledge of set-theoretic defining polynomials for
$V(T_3)$ is sufficient to give set-theoretic defining polynomials
for $V(T)$ for any binary tree $T$. Thus while one might naively
view the case of $V(T_3)$ as the simplest, in fact it is the only
remaining barrier to the determination of polynomials defining the
binary $n$-taxon variety, for any $n$. In the cases $\kappa=2,3$
where such defining polynomials are known, we thus obtain the
following.

\begin{cor}
For $\kappa=2$ or 3, and any binary tree $T$, explicit polynomials
set-theoretically defining V(T) can be given.
\end{cor}

Note that for $\kappa=2$ a stronger result is provided by Theorem
\ref{thm:PS}, to be proved in Section \ref{sec:k2id}.

\medskip

For the remainder of this section let $V_{Flat}(T)$ denote the zero
set of  $\mathcal F(T)$.  Our proof of Theorem \ref{thm:kany} will
follow several lemmas. The first is an analog for $V_{Flat}(T)$ of
Lemma \ref{lem:closst}.

\begin{lem}\label{lem:fvalg} Let $T'$ and $T''$ be $n$-taxon and
$m$-taxon trees, with $T=T'\st T''$.  If $Q\in CV_{Flat}(T')$ and $R\in
CV_{Flat}(T'')$, then $Q\st R \in CV_{Flat}(T)$.  \end{lem}

\begin{pf} Consider any internal node $v$ of $T$, which we may assume
arises from an internal node of $T'$. We assume $v$ is trivalent;
straight-forward modifications to our argument give the general
case.

Flattening $Q$ at $v$, the resulting tensor lies on $CV(\kappa;
\kappa^{n_1}, \kappa^{n_2}, \kappa^{n_3})$, with $n_3=n-n_1-n_2$,
where we assume taxon $a_n$ of $T'$ (where taxon $b_1$ of $T''$ is
to be joined) is included in the last index of the flattening.  Then
the flattening of $Q\st R$ at $v$ is obtained from the flattening of
$Q$ at $v$ by an action in the third index by a matrix $R'$ whose
entries are determined by those of $R$. By Lemma \ref{lem:extedge}
the flattening of $Q\st R$ at $v$ lies in $CV(\kappa;
\kappa^{n_1},\kappa^{n_2},\kappa^{n_3+m-2})$.

Thus $Q\st R\in CV_{Flat}(T)$.  \qed\end{pf}

We also need a converse to this lemma.

\begin{lem} \label{lem:fvdecomp} Let $T'$ and $T''$ be $n$-taxon and
$m$-taxon trees, with $T=T'\st T''$.  Then if $P\in CV_{Flat}(T)$, there
exist $Q\in CV_{Flat}(T')$ and $R\in CV_{Flat}(T'')$ with $P=Q \st R$.
\end{lem} \begin{pf}

Let $e$ be the edge of $T$ formed by conjoining edges of $T'$ and $T''$.
Since any $P\in CV_{Flat}(T)$ satisfies the edge invariants for $e$, we
may flatten it on $e$ to obtain a $\kappa^{n-1}\times \kappa^{m-1}$
matrix of rank $l\le \kappa$, and write $$P=Q \st R,$$ where $Q$ and $R$
are $n$- and $m$-dimensional tensors, respectively, with all indices
running through $[\kappa]$.  We may further assume the non-zero
$Q_k=Q(\cdot,\cdots,\cdot,k)$ are linearly independent, as are the
non-zero $R_k=R(k,\cdot,\cdots,\cdot)$, and that $Q_k,R_k$ are non-zero
only for $k=1,\dots, l\le\kappa$.

We next show $Q\in CV_{Flat}(T')$.  First observe that since the
non-zero $R_k$ are independent, if we write them as row vectors,
there is a $\kappa^{m-1}\times \kappa$ matrix $A$ so that
$R_kA=\mathbf e_k$ for all $k\le l$.  Now supposing the taxa of $T'$
and $T''$ are $a_1,\dots, a_n$ and $b_1,\dots, b_m$, respectively,
flatten $P$ according to the partition $\left \{ \{a_1\},\dots,
\{a_{n-1}\},\{b_{2},\dots, b_m\} \right \}$ to an $n$-dimensional
$\kappa\times\cdots \times\kappa\times\kappa^{m-1}$ tensor $F$.
Letting $R'$ denote the $\kappa\times\kappa^{m-1}$ flattened form of
$R$ with rows $R_k$, we have $F=Q*R'$. Thus $F*A=Q*R'*A=Q$.  (Note
that $A$ does not act in a single index of $P$ here, but does act in
a single index of the flattening $F$.) It is now straightforward to
see that any flattening of $Q$ at an internal vertex of $T'$ is
obtained from a flattening of $P$ at a vertex of $T$, followed by an
action in one of the resulting indices of a matrix determined by
$A$. Thus by the definition of $\mathcal F(T)$, $Q$ will satisfy all
polynomials in $\mathcal F(T')$.

Similarly, $R\in CV_{Flat}(T'')$.  \qed\end{pf}

\begin{pf}[Proof of Theorem \ref{thm:kany}] We already know that
$V_{Flat}(T) \supseteq V(T)$.

The proof that $V_{Flat}(T)= V(T)$ proceeds by induction on the number
$n$ of taxa.  The cases of $n=2,3$ hold by the definition of $\mathcal
F(T)$.

For simplicity, we first consider a binary tree $T=T_n$,with $n\ge
4$ taxa. Picking a cherry of $T$, let $T_{n-1}$ and $T_3$ be such
that $T=T_{n-1}\st T_3$. Suppose $P\in CV_{Flat}(T)$. By Lemma
\ref{lem:fvdecomp}, we have $P=Q\st R$, for $Q\in
CV_{Flat}(T_{n-1})$ and $R\in CV_{Flat}(T_3)$. This, in combination
with Lemma \ref{lem:fvalg}, means the map
$$\mu: CV_{Flat}(T_{n-1})\times CV_{Flat}(T_3) \to CV_{Flat}(T_n)$$
defined by $(Q,R)\mapsto Q\st R$ is surjective.

Denote the parameterizations of the cones over the phylogenetic
varieties for $T_k$ by $\psi_k:U_k\to\C^{L_k}$.  With the map
$\alpha:U_{n-1}\times U_3 \to U_n$ defined by
$\alpha(u_{n-1},u_3)=u_{n-1}\st u_3,$ the diagram $$ \begin{CD}
U_{n-1}\times U_3 @>{\psi_{n-1} \times \psi_3}>>
CV_{Flat}(T_{n-1})\times CV_{Flat}(T_3)\\ @V\alpha VV @V\mu VV \\ U_n
@>\psi_n>> CV_{Flat}(T) \end{CD} $$ commutes, by Lemma
\ref{lem:starparam}.

Now $\alpha$ and $\mu$ are surjective, and by the inductive
hypothesis the image of $\psi_{n-1}\times\psi_3$ is dense in
$CV_{Flat}(T_{n-1})\times CV_{Flat}(T_3)$, so the image of $\psi_n$
is dense in $CV_{Flat}(T)$. Thus $V_{Flat}(T)= V(T)$.

If $T$ is not binary, the above argument may be modified by
replacing the decomposition $T=T_{n-1}*T_3$ by $T=T_{n-k+1}*T_{k+1}$
where $T_{k+1}$ is a star tree with $k+1$ leaves and $T_{n-k+1}$ has
$n-k+1$ leaves, if necessary. \qed\end{pf}

Theorem \ref{thm:kany} and the preceding lemmas yield the following
strengthening of Lemma \ref{lem:closst}.

\begin{cor}\label{cor:pvalg} If $T=T'\st T''$, then $CV(T)=CV(T')\st
CV(T'')$.  \end{cor}

\section{The phylogenetic ideal: Binary $T$ and $\kappa=2$. }\label{sec:k2id}

We now prove Theorem \ref{thm:PS}, and thus assume $T$ is a binary
tree and $\kappa=2$.

Our arguments will use in several ways the fact (see Section
\ref{sec:startrees}) that for $\kappa=2$ the variety $V(T_3)$ fills
its ambient space: $V(T_3)=\PP^{7}.$ Note, however, that for
$\kappa>2$, $V(T_3)\subsetneq\PP^{\kappa^3-1}$, and so the approach
here cannot be successfully modified in a simple way.

The first use of this special fact is to note that for our chosen
$\kappa$, $V(2;2,2,2)=V(T_3)=\PP^7$ means the set $\mathcal F$
defining $V(T_3)$ is $\{0\}$. Thus the set $\mathcal F(T_n)$ of the
set-theoretic result Theorem \ref{thm:kany} is the set of edge
invariants.  While our goal is to show $\mathcal F(T_n)$ generates
the full ideal vanishing on $V(T_n)$, we will not, in fact, appeal
to Theorem \ref{thm:kany} to do so.

The second use of $V(T_3)=\PP^{7}$ is more subtle.  Recall that
regardless of $\kappa$, there are actions of $GL(\kappa,\C)$ on
$V(T_n)$ in each index. However, in the case $\kappa=2$, the special
nature of $V(T_3)$ gives us actions of $GL(4,\C)$ on $V(T_n)$ via
the cherries of $T_n$.  This is really the key point in our
argument, as it underlies the application of Lemma \ref{lem:reparg}.
Nonetheless, this action is in some respect an `unnatural'
consequence of $\kappa=2$.  The following lemma provides a more
careful statement of the special structure we use.

\begin{lem} \label{lem:24st}
Let $T_n$ denote a binary $n$-taxon tree, labeled so that taxa
$a_{n-1}$ and $a_n$ form a cherry. Write $T_n=T_{n-1}\st T_3$, where
$a_{n-1}, a_n$ are taxa on $T_3$. Let $e$ denote the edge of $T_n$
formed from conjoining edge $\tilde e$ of $T_{n-1}$ and the
appropriate edge of $T_3$.  View points in $CV(T_n)$ and
$CV(T_{n-1})$ as $2^{n-2}\times 4$ and $2^{n-2}\times 2$ matrices by
flattening them on the edges $e$ and $\tilde e$, respectively.  Then
$$CV(T_n)=\overline{CV(T_{n-1})\st M(2,4,\C)}$$ and
$$CV(T_{n-1})=CV(T_{n-1})*M(2,2,\C).$$ \end{lem}

\begin{pf} The first claim is simply Lemma \ref{lem:closst} applied
to $T_{n-1}$ and $T_3$, combined with the observation that
$CV(T_3)=\C^8$ flattens to give $M(2,4,\C)$.  (Note that by Corollary
\ref{cor:pvalg}, we could also remove the closure symbol here.)

For the second claim, apply the same argument to $T_{n-1}$ and
$T_2$, observing that $CV(T_2)=M(2,2,\C)$.  \qed\end{pf}

\begin{pf}[Proof of Theorem \ref{thm:PS}] We proceed by induction on
the number $n$ of taxa for $T_n$, with the cases of $n=2,3$ known.

Let $\aideal=\aideal_T$ denote the ideal vanishing on $CV(T)$, and
$\bid$ the ideal generated by $\mathcal F_{edge}(T)$.  That
$\aideal\supseteq \bid$ has been discussed already; we must show the
opposite inclusion.

With $T_n=T$, choose a cherry so that $T_n=T_{n-1}\st T_3$, with
notation as in Lemma \ref{lem:24st}. By that lemma, we may apply
Lemma \ref{lem:reparg} with $V_1=CV(T_{n-1})$ and $V_2=CV(T_n)$.  We
thus find $\aideal$ is generated by the $3\times 3$ minors of the
edge flattening $Flat_e(P)$ on the conjoined edge $e$ of an
$n$-dimensional tensor of indeterminates $P$, together with all
polynomials of the form $g(P)=f(Flat_e(P)*B)$ where $f(Q)$ vanishes
on $CV(T_{n-1})$, $Q$ is a $(n-1)$-dimensional tensor of
indeterminates, and $B\in M(4,2,\C)$.

Now, by induction, the ideal of such $f$ is generated by $3\times 3$
minors of $Flat_{e'}(Q)$ as $e'$ ranges through edges of $T_{n-1}$.
Consider one such minor, say $f_0$, obtained from the flattening on
an edge $e_0$ of $T_{n-1}$. We may assume $e_0\ne \tilde e$,
since otherwise there are no $3\times3$ minors.
It will be enough to show
$f_0(Flat_e(P)*B)\in \bid$.

We claim that $f_0(Flat_e(P)*B)$ vanishes on all $P$ that have rank at
most 2 when flattened on the edge $e_0$ in $T_n$. For such a $P$, since
$Flat_{e_0}(P)$ is $2^m\times 2^{n-m}$, there is an expression
$P=P_1*P_2$, where $P_1$ is an $(m+1)$-dimensional $2\times\cdots\times
2$ tensor, and $P_2$ an $(n-m+1)$-dimensional $2\times\cdots\times 2$
tensor.  Then writing $P$ and $P_2$ as $2\times \cdots \times 2\times 4$
tensors by flattening to combine the taxa $a_{n-1},a_n$, we have
$P*B=P_1*(P_2*B)$. This shows $P*B$ also has rank at most 2 when
flattened on $e_0$, and so $f_0$ vanishes on it, as claimed.

But since $f_0(Flat_e(P)*B)$ vanishes on all $P$ of rank at most 2
when flattened on $e_0$, it is contained in the ideal generated by
$3\times 3$ minors of flattenings on $e_0$. Thus it is in $\bid$.
\qed\end{pf}

\bibliographystyle{plain}
\bibliography{phylo}

\end{document}